\documentclass{article}
\usepackage{graphicx}
\usepackage{natbib}
\def\G{{\cal G}}

\begin{document}

%%%%%%%%%%%%%%%%%%%%%%%%%%%%%%%%%%%%%%%%%%%%%%%%%%%%%%%%%%%%%%%%%%%%%%%%%%%
%\markboth{K. Judd}{Non-probabilistic odds and forecasting with imperfect models}
\title{Non-probabilistic odds and forecasting with imperfect models}

%%%%%%%%%%%%%%%%%%%%%%%%%%%%%%%%%%%%%%%%%%%%%%%%%%%%%%%%%%%%%%%%%%%%%%%%%%%
\author{Kevin Judd}
%\address{School of Mathematics and Statistics, M019,\\
%  University of Western Australia, Perth, 6009, Australia. \email{Kevin.Judd@uwa.edu.au}}

\maketitle

%%%%%%%%%%%%%%%%%%%%%%%%%%%%%%%%%%%%%%%%%%%%%%%%%%%%%%%%%%%%%%%%%%%%%%%%%%%
\begin{abstract}
  Probability forecasts are intended to account for the uncertainties
  inherent in forecasting. It is suggested that from an end-user's point
  of view probability is not necessarily sufficient to reflect
  uncertainties that are not simply the result of complexity or
  randomness, for example, probability forecasts may not adequately
  account for uncertainties due to model error. It is suggested that an
  alternative \emph{forecast product} is to issue non-probabilistic odds
  forecasts, which may be as useful to end-users, and give a less
  distorted account of the uncertainties of a forecast.  Our analysis of
  odds forecasts derives from game theory using the principle that if
  forecasters truly believe their forecasts, then they should take bets at
  the odds they offer and not expect to be bankrupted.  Despite this game
  theoretic approach, it is not a market or economic evaluation; it is
  intended to be a scientific evaluation.  Illustrative examples are given
  of the calculation of odds forecasts and their application to
  investment, loss mitigation and ensemble weather forecasting.
\end{abstract}

%\begin{keywords}
%  non-probabilistic odds; probability forecasting
%\end{keywords}

%%%%%%%%%%%%%%%%%%%%%%%%%%%%%%%%%%%%%%%%%%%%%%%%%%%%%%%%%%%%%%%%%%%%%%%%%%%

\section{Introduction}
% Our purpose here is to consider the quantification of uncertainty in
% forecasting. Our main concern here is whether probability necessarily
% captures all aspects of uncertainty in forecasting, and whether
% non-probabilitic odds may provide a useful alternative.
%
% ``Put your money where your mouth is'' is a powerful taunt: if you really
% believe what you say, then you should be prepared to bet your hard-earned
% wealth on it. A wager is not an uncivilised way of resolving a dispute.
% When antagonists negotiate a value and odds for the wager, they quantify
% their certainty.

Our concern here is the quantification of uncertainty in forecasting.
Suppose our task is to forecast the possibility of overnight temperatures
at a particular location falling below freezing. This is a task of
significant economic importance, for example, for the salting and gritting
of roads, and for preparations to mitigate of frost damage in
horticulture. Forecasts of this type use atmospheric observations and
computer models of the physical processes of the atmosphere. Uncertainty
arises in these forecasts from the complexity of the weather, from the
sparsity of observations, from the simplifications and inadequacies of the
computer models, and so on. The suggestion we make here is that in some
situations, such as using forecasts to mitigate losses, probability may
not be the best means to quantify the uncertainty of the forecast. We
suggest that non-probabilistic odds may provide a useful alternative to
some forecast users.

A probability forecast quantifies uncertainty by assigning a probability
to the occurrence of an event. If the process that determines the outcome
of the event is intrinsically random, then there can be a \emph{correct}
or \emph{optimal} assignment to the probability. The forecaster, however,
may be uncertain about the correct probability value to assign, and this
leads to consideration of the influence of other uncertainties, which we
discuss shortly. A well established means for quantifying the uncertainty
of probability forecasting systems are reliability or skill scores.
\cite{Brier:1950}, Murphy, Winkler \citep{Murphy-Winkler:1977,
  Murphy-Winkler:1987, Murphy:1993, Winkler:1994}, and many others since,
have considered using scores to assess the skill of a probability
forecast.  Skill scores are closely related to issues of calibration of
probability forecasts~\citep{Foster-Vohra:1998, Palmer:2000,
  Roulston-Smith:2002, Smith:1995}. An alternative assessment of skill is
the economic value of the forecast~\citep{Granger-Pesaran:2000}.  Economic
value is closely related to the idea of wagers. One can imagine a market
of forecasters who take bets on the outcomes of events they forecast; the
forecasters that profit the most, or at least avoid bankruptcy, are
considered the better forecasters. In a competitive market only those
forecasters who know the true probabilities of events will survive in the
long term~\citep{Shafer-Vovk:2001}.  In the short term, however, merely
lucky forecasters can survive, and even excel~\citep{Johnstone:2007}.
Consequently, for a forecaster to be a top performer they are forced to
act like a bookmaker or marketeer rather than a scientist.  Scientists aim
to learn the true probabilities of events, whereas bookmakers and
marketeers respond to the opportunities of a market~\citep{Levitt:2004}.

% (It will be necessary later to isolate forecaster's wagers from market
% information to avoid scientists being forced into being bookmakers or
% marketeers.)

A probability forecast and a skill score %(or economic value)
together give a more complete picture of the uncertainty of a forecasted
event. We will argue here, however, that probability forecasts can give
end-users a misleading picture of uncertainty. We suggest a method of
avoiding this problem is to issue odds forecasts, which arise naturally as
wagers. In the next subsections we state our distinction between
\emph{odds} and \emph{probability}, and describe briefly how odds can
quantify multiple aspects of uncertainty. Section~\ref{sec:maths}
introduces a mathematical formalism for the process of issuing odds
forecasts, including a brief review of the necessary concepts and
techniques of game theory. Most importantly we show how the computation of
odds can be framed as a simple optimisation problem.  In
section~\ref{sec:models} we compute odds in four basic forecasting
situations and compare these computed odds with the corresponding
probabilities. In section~\ref{sec:investment} we describe how odds can be
employed in investment and loss mitigation.  Finally, in
section~\ref{sec:ensembles} we use operational ensemble weather forecasts
and station data to present a concrete example of issuing odds forecasts.

\subsection{Odds and probability}\label{sec:odds}
% Since the eighteenth century, as a consequence of work of Jacob Bernoulli,
% De Moivre, Simpson, Bayes and Laplace~\citep[see][]{Stigler:1986},
% probability has been the quantity used by scientists to express
% uncertainty. In keeping with this, many forecasters want to make
% \emph{probability forecasts}~\citep[for example:][]{Murphy-Winkler:1977, Murphy-Winkler:1987, Murphy:1993,
%   Winkler:1994, Pan-Dool:1998}.  Probability forecasts can be used for
% decision making~\citep{Stewart:1997, Mylne:2002} and forecast
% evaluation~\citep{Smith:1995, Palmer:2000,
%   Granger-Pesaran:2000, Roulston-Smith:2002}.
%
%We argue here that a probability forecast is not always the
%most appropriate forecast and propose \emph{non-probabilistic odds} as an
%alternative. Part of the reason for this view is that there are different
%kinds of uncertainty.

% Odds have been used by bettors, bookmakers and casinos for centuries and
% are still used today.  One traditional way to express odds is by their
% \emph{pay-out}; odds of ``$A$ to~$B$'', means a pay-out of $A$~dollars on
% a bet of $B$~dollars. Odds can be interpret as an expression of
% uncertainty, using the quantity~$q=B/(A+B)$.  Often odds and probability
% are treated as though they were equivalent, but they are equivalent only
% in specific circumstances.
%
Given a complete set of $m$ mutually exclusive events, that is, one and
only one of the events will occur, we define \emph{odds} to mean an
assignment of real numbers $q_i\geq0$, $i=1,\dots,m$, to the events. If
$\sum_{i=1}^mq_i=1$, then the odds are \emph{probabilistic odds}, and the
$q_i$ are probabilities.
%In our terminology, the
%pay-out, ``$(1/q)-1$ to $1$'' or ``$1-q$ to~$q$'', is just an alternative
%way of expressing odds.  
Casinos and bookmakers assign odds so that $\sum_{i=1}^mq_i>1$, which has
the consequence that, on average, they should profit from the bettors. The
excess over one of the sum is sometimes termed ``juice'', ``take'',
vigour''.

%Betting plays an important role in our analysis. 
We argue that if a forecaster truly believes their forecasts, then they
should take bets at the odds they offer and not expect to be bankrupted.
This can always be achieved by a sufficiently large excess, however, the
scientist's goal ought to be to avoid bankruptcy with the smallest excess.

Taking wagers while avoiding bankruptcy is a powerful principle. Extending
early work of Ville it has been used by
%\citet{Foster-Vohra:1998}, 
%\citet{Schervish:1989},
\citet{Foster-Vovk:1999}, 
\citet{Skouras-Dawid:1999},
%\citet{Granger-Pesaran:2000}, 
\citet{Shafer-Vovk:2001}, 
\citet{Dawid:2004}
and others.
%(Notes on which this
%paper is based have been circulated for sometime and the author became
%aware of this other work rather belatedly.)  
\citet{Shafer-Vovk:2001} provides a beautiful development that derives
probability theory itself from betting principles. The ideas developed in
this paper share a conceptual and structural formulation with
\citet{Shafer-Vovk:2001}, but there are significant differences.  Briefly,
the differences arise because in our opinion non-probabilistic odds are
relevant to the immediate and short-term consequences of uncertainty,
whereas probabilities are generally more relevant to long-term and asymptotic
uncertainty.
%
% In these notes we first investigate how a forecaster should assign odds to
% events so that even an adversary who knew the true probabilities (or
% asymptotic frequencies) of events would not bankrupt the forecaster. We
% demonstrate the calculation of such odds in a number of different
% scenarios. Finally we discuss how such non-probabilistic odds forecasts
% can be used by an investor to mitigate their losses. It is significant to
% note that investors work directly with the odds, there is no need to
% normalise the odds in an attempt to interpret them as estimates of
% probabilities.

Bid-ask spreads in financial markets can be interpreted as
non-probabilistic odds~\citep{Levitt:2004, Johnstone:2007}, but our
development of odds differs because we want our forecaster to behave as a
scientist and avoid market pressures influencing the odds.

\subsection{Levels of uncertainty}

Uncertainty arises from many sources, not just randomness.
% For example, \emph{first order} uncertainty is
% represented by the probabilities assigned to events, where as, \emph{second
%   order} uncertainty represents uncertainties about the probabilities that
% have been assigned.  One might consider higher order uncertainties, which
% might represent uncertainties about the appropriate probabilistic model
% of the events.  In these notes we suggest that non-probabilistic odds can
% conveniently represent both first and second order uncertainty under the
% assumption that there are not higher order uncertainties.
To appreciate this consider Jacob Bernoulli's foundational example of
drawing balls, with replacement, from a urn containing a number of red and
black balls. Suppose the fraction~$\theta$ of red balls in the urn is
known, and the process of extracting a ball is sufficiently complex to
appear \emph{uniformly random}, that is, on any selection each ball in the
urn is equally likely to be drawn.  In this situation the uncertainty is
entirely due to the ``random'' process of selection, and the uncertainty
is adequately described by the probability of drawing a red ball, which is
$\theta$. We will refer to this situation as having \emph{first order}
uncertainty.

Knowing $\theta$ allows one to make a probability forecast, that is, the
probability of the event ``drawing a red ball'' is~$\theta$. This is a
simple \emph{forecast model}. Indeed, under our assumptions, this is a
\emph{perfect model}, because the model exactly represents the process,
and there is no better model.  In a perfect model there is only first
order uncertainty.

If the number of red and black balls in the urn is unknown, then a
value~$\psi$ for the fraction of red balls could be inferred from the
observed fraction of red balls seen in a number of draws. Using $\psi$ to
forecast the probability of a red ball is an \emph{imperfect model}, there
is both first order uncertainty from the randomness
%(or complexity) 
of the process being forecast, and second order
uncertainty due to the value of $\psi$ used.

It is possible to make a further distinction between a \emph{perfect model
  class} and an \emph{imperfect model class}. For the urn a value of
$0\le{}\psi\le{}1$ defines a class of forecast models. If the assumption
of uniform random selection of balls is valid, then this is a perfect
model class, because the value $\psi=\theta$ provides a perfect model. If
the random selection assumption does not hold, because perhaps the balls
are not thoroughly stirred before each selection, then no value of $\psi$
provides a perfect model; a perfect model would have to take into account
the conditional randomness, or even non-random effects, of mixing the urn
after each replacement.  Uncertainty about the correct model class is at
least third order.  
%Much of statistical theory implicitly assumes a
%perfect model class.

\subsection{At odds with probability}
One goal of forecasters is to issue accountable, or reliable, or well
calibrated, probability forecasts~\citep{Foster-Vohra:1998, Smith:1995,
  Palmer:2000}, that is, if it is forecast that an event~$E$ will occur
with probability~$p$, then the observed fraction of events~$E$ out of all
events asymptotically approaches $p$. It is vital to recognise that an
imperfect model will almost certainly not provide an accountable
probability forecast. At best forecasts are only accountable asymptotically
in a perfect model class, but even this is a delicate
problem~\citep{Oakes:1985, Foster-Vohra:1998}.

Suppose that for an urn game a forecaster a provides odds for outcomes of
draws, and takes bets at these odds. If the forecaster knew $\theta$, and
the uniform random selection assumption held, then probabilistic odds
could be issued where a one unit bet on a red ball pays $1-\theta$ and
betting a black ball pays $\theta$.  With these probabilistic odds the
forecaster does not expect their wealth to increase or decrease on
average, no matter how skillful the bettors. If the forecaster had only an
estimate $\psi$, then providing probabilistic odds, where a bet on a red
ball pays $1-\psi$ and a bet on a black ball pays $\psi$, would be unwise.
Any bettor who knew $\theta$ would almost surely bankrupt the forecaster
by always betting a fraction~$\theta$ of their current wealth on the red
ball, and betting the fraction~$1-\theta$ of their current wealth on the
black ball~\citep{Kelly:betting}. A similar result is true when the bettor
just has a better estimate of~$\theta$.
%(We will develop this result later.)
The fact that $\psi=\theta$ is the only value of $\psi$ that does not lead
asymptotically to certain bankruptcy can be used to define the concept of
probability~\citep{Shafer-Vovk:2001}.

The essential problem with using imperfect models to provide probability
forecasts is the higher order uncertainties, like model error, distort the
odds. In particular, we will show that model error often results in under
estimating the chance of events that have low probability, a so-called
base-rate effect. We argue that by using non-probabilistic odds forecasts,
rather than probability forecasts, a forecaster can provide a less
distorted forecast that takes into account model error.  (Here we
only show how to do this for particular second order uncertainties.)

We argue that if a forecaster aims to provide odds that are as close to
probabilistic as possible, subject to known information and the model
class used, then the amount by which the odds exceed one is a measure of
how certain the forecaster is of their model. The excess is a measure of
second order uncertainty.  (Here we only demonstrate how to do
this for a perfect model class with particular utility functions.)

%%%%%%%%%%%%%%%%%%%%%%%%%%%%%%%%%%%%%%%%%%%%%%%%%%%%%%%%%%%%%%%%%%%%%%%%%%%
\section{Computation of odds forecasts}\label{sec:maths}
The theoretical framework in which we analyse odds forecasting is game
theory~\citep{game-theory, Shafer-Vovk:2001}. There are three
``players'' in our game: the Forecaster, the Client, and Nature.  The
Client determines the structure of the game.  Imagine that the Client
approaches the Forecaster and requests odds for a complete set of $m$
mutually exclusive events $E_i$, $i=1,\dots,m$. Which event actually
occurs is determined by Nature, who is ignorant and indifferent to the
Client-Forecaster negotiations. For our purposes we may assume that
Nature's selection of an event is random, or of such complexity that the
selection is assumed random by the Client and Forecaster. (By making this
assumption about Nature we avoid issues of third order uncertainties.)

When the Client approaches the Forecaster, the Client does not specify
interest in any particular event within the set of mutually exclusive
events. The Forecaster is required to supply odds for all events, based on
past observations of Nature and the Forecaster's model.  The Forecaster is
not a bookmaker, they are a scientist. The objective of the
scientist Forecaster is to provide odds as close to probabilistic as
possible given the available information and their model class, because by
doing so they are aiming to obtain the best forecast model.

It may be tempting to imagine many clients and forecasters competing in a
market to determine the best forecast model, but this has the propensity
for forecasters to adopt the strategies of marketeers and bookmakers.
To avoid forecasters acting this way we will isolate the Forecaster from
market information: there will be one Client and the Forecaster is not
allowed to know how the Client bets, only the choices of Nature. If the
Forecaster knew the pattern of the Client's bets, then this constitutes an
additional information stream~\citep{Kelly:betting}, and consequently can
be used to improve the Forecaster's performance against the Client, or
against Nature if the Client is a more skillful forecaster than the
Forecaster. Using this additional information stream allows marketeering
and bookmaking, rather than scientific forecasting.  Denying the
Forecaster knowledge of the Client's bets forces the Forecaster to rely on
available observations and modelling skill alone. Shortly we will see,
however, that the Forecaster will need to know a little about the Client's
betting.

\subsection{The Game}
%The following is the nature of the game. %~\citep{RollingStones:Sympathy}.
We have a three person game. Nature plays indifferently with no aim.  The
Client aims to accumulate winnings. The Forecaster aims to set as close to
probabilistic odds as possible, without the Client bankrupting the
Forecaster. If there were just two events, $E$ and its complement~$E'$,
then the game can be represented by a game matrix~$\G$.
\begin{equation}
  \label{eq:game2x2}
  \G\colon\begin{tabular}{cc}
      \raisebox{-0.5cm}{client} & 
      \begin{tabular}{c|cc|}
        \multicolumn{1}{c}{} & \multicolumn{2}{c}{nature}\\
        & $E$ & $E'$ \\
        \hline
        $E$ & $P$ & $-1$ \\
        $E'$ & $-1$ & $P'$ \\
        \hline
      \end{tabular}
    \end{tabular}
\end{equation}
The game matrix represents the pay-out to the Client for a one unit bet on
an event given Nature's outcome. The odds set by the Forecaster in this
case are $P$ to~$1$ for event~$E$, and $P'$ to~$1$ for event~$E'$.

% Operationally the Client is playing a perfect information, zero-sum game
% against Nature, using a pay-out matrix specified by the Forecaster, but
% conceptually, the Client is playiith:2002}.
%ng a zero-sum game against the
% Forecaster. On the other hand, conceptually, the Forecaster is trying to
% provide as probabilistic odds as possible, which translates operationally
% into setting odds (without knowing how the Client plays) so that on
% average the Client does not accumulate winnings, while setting the odds so
% that the Client is still motivated to play.  The game is not a perfect
% information game; the Client has the information advantage, but the
% Forecaster has the advantage of being able to set non-probabilistic odds.

There are several variants of this game according to the rules that govern
the Client's bets; the variants can influence how the Forecaster sets the
odds. One rule that can be introduced is the Client's bets are a fixed
size and a negligible fraction of the Client's (and Forecaster's) total
wealth, effectively, infinitesimal bets.  Alternatively, the Client can
bet a substantial fraction of their wealth.  Other rules that might be
introduced govern whether the Client is forced to make a bet, regardless
of the odds, or whether the Client can split bets over several events, or
whether there is a minimum size of a bet on any event. Rules governing the
sizes of bets can influence the Client's utility function of wealth, and
consequently influence the pattern of bets. With variable sized bets the
Client can choose to maximise the growth of wealth, a logarithmic utility
function, which requires distributing bets over several events. Forced
bets of a fixed (infinitesimal) size essentially implies a linear utility
function of wealth, and results in betting only on the event that the
Client believes gives the maximum pay-out.

Since the Forecaster does not know the Client's individual bets, the
Forecaster must at least know the Client's utility function of wealth, or
the rules governing how the Client can bet, which effectively force a
utility function onto the Client.

The ultimate challenge for the scientist Forecaster is to compete against
a Client who knows the true probabilities of Nature, and yet do so without
being bankrupted.

\subsection{Forced bets of fixed (infinitesimal) size; linear utility}\label{sec:forced}
 Much of the analysis in this section is text book game
theory~\citep{game-theory}.
%Should clarification be required, the simplest
%case of a single event and its complement is worked in the appendix. 
We first consider the zero-sum game between the Client and Nature and
determine the Client's optimal strategy for fixed odds. We then derive an
optimal odds assignment of the Forecaster for an arbitrary probability
model under the assumptions of a perfect model class.

\subsubsection{Optimal Client strategy}
For a complete set of $m$ mutually exclusive events $E_i$, $i=1,\dots,m$,
the $m\times{}m$ game matrix~$\G$ of the zero-sum game between the Client
and Nature (represented in $2\times{}2$ case in~(\ref{eq:game2x2})) has
$\G_{ij}=(\delta_{ij}/q_i) - 1$, where $\delta_{ij}=1$ if $i=j$ and zero
otherwise, and $q_i$ is the odds for event~$E_i$, or equivalently, the
odds on event~$E_i$ is ``$P_i$ to~$1$'', where $q_i = 1/(P_i+1)$. Assume
that Nature plays a strategy where event~$E_i$ is chosen randomly with
probability~$\pi_i$. Supposing that the Client plays a strategy where
event~$E_i$ is chosen with probability~$p_i$, then the average pay-out to
the Client is 
\begin{equation}
  \label{eq:payout}
  \left(\sum_{i=1}^m\frac{\pi_ip_i}{q_i}\right) - 1.  
\end{equation}
Here $\sum_{i=1}^m\pi_i=1$, $\sum_{i=1}^mp_i=1$, but
$\sum_{i=1}^mq_i\ge{}1$, with equality if and only if the odds are
probabilistic.

The \emph{minimax} strategy for the Client is to select event~$E_i$ with
probability $p^\star_i$, and has average pay-out $V^\star$, where
\begin{displaymath}
  \label{eq:minimax}
  p^\star_i=\frac{q_i}{\sum_iq_i}, \quad 
  V^\star=\frac{1}{\sum_iq_i}-1.
\end{displaymath}
By using the minimax strategy the Client is guaranteed an average pay-out
of at least $V^\star$ regardless of the~$\pi_i$. When the odds are
probabilistic, then the minimax strategy for the Client has $p_i=q_i$ and
$V^\star=0$.  It follows that the optimal Client strategy is to choose
$p_i=1$, when $i$ is the index such that $\pi_i/q_i$ is maximal, and
$p_i=0$ otherwise.  The average pay-out to the Client is then
$\max_i\{\pi_i/q_i\}-1$.

\subsubsection{Optimal Forecaster strategies}
The forecaster's aim is to make the odds as probabilistic as possible
without the client's winnings accumulating without bound. The forecaster
has to allow for the possibility that the client is a better forecaster,
in the worst case, the Client knows Nature's
probabilities~$\pi=(\pi_1,\dots,\pi_m)$. There is no single optimal
strategy to achieve the goals we have set for the forecaster, so some
additional guidance is necessary.

%Making the odds as probabilistic as possible can be taken to mean minimising
%$Q=\sum_{j=1}^mq_j$.  
Since the game matrix~$\G$ is symmetric, it can be seen from
equation~\ref{eq:payout} that the \emph{ideal} strategy for the forecaster
is to set the odds so that $q_i/\sum_{j=1}^mq_j=\pi_i$, in which case the
closest to probabilistic odds are the probabilistic odds $q_i=\pi_i$. This
is, of course, a meaningless solution, because the $\pi_i$ are unknown.
Simply taking an estimate~$\hat\pi$ of Nature's probabilities~$\pi$, and
using these as odds, $q_i=\hat\pi_i$, is unwise, because at least one
$\pi_i/\hat\pi_i>1$, which will be exploited by the more informed Client,
who can obtain a better estimate of~$\pi_i$.

If the forecaster accepts their forecast model is imperfect, then they
will acknowledge there are many likely values of~$\pi$. However, if the
forecaster assumes their model class is perfect, then given data~$D$ they
can assert that under their model $\Pr(\pi\mid{}D)\propto\Pr(D\mid{}\pi)\Pr(\pi)$,
where $\Pr(\pi)$ represents prior knowledge about possible values of
$\pi$.  From this the forecaster can compute the average loss (pay-out to
the client) $V(\pi,q)$ for given~$\pi$ and odds $q=(q_1,\dots,q_m)$,
\begin{equation}
  \label{eq:V}
  V(\pi,q) = \max_i\left\{\frac{\pi_i}{q_i}\right\}-1
\end{equation}
The expected loss $E(V\mid{}D,q)$, given the data~$D$ and fixed odds~$q$, is
\begin{equation}
  \label{eq:EV}
  E(V\mid{}D,q) = \frac{\int_SV(\pi,q)\Pr(D\mid{}\pi)\Pr(\pi)\,d\pi}{\int_S\Pr(D\mid{}\pi)\Pr(\pi)\,d\pi},
\end{equation}
where $S$ is the simplex $\sum_i\pi_i=1$.  A possible optimal strategy for
the forecaster is
\begin{equation}
  \label{eq:fopt}
  \min \sum_iq_i \quad\mbox{subject to}\quad 
  \int_SV(\pi,q)\Pr(D\mid{}\pi)\Pr(\pi)\,d\pi=0.
\end{equation}
The minimisation attempts to ensure the odds~$q$ are as probabilistic as
possible.  The constraint $E(V\mid{}D,q)=0$ attempts to ensure that whatever
the true~$\pi$, the odds are such that the client's average pay-out is
zero given the assumed distribution of~$\pi$. The construction of
equation~(\ref{eq:fopt}) implies that the client's wealth is a
\emph{martingale}~\citep{Doob:stochastic-processes,Williams:martingales},
that is, if $W_n$ is the wealth of the client after $n$ plays of the game,
then $E(W_n\mid{}D,q) = W_{n-1}+E(V\mid{}D,q)= W_{n-1}$.

% The forecaster does not guarantee the client's winnings will not grow
% without bound, only that they are not expected to. The strategy given by
% equation~(\ref{eq:fopt}) is not claimed to be the best strategy for the
% forecaster, but it seems not to be an unreasonable strategy, and is
% certainly better than a probability forecast. In next section we
% investigate the odds it provides for a frequency and a Gaussian model.

Observe that the constraint $E(V\mid{}D,q)=0$ is not, in principle, hard to
compute. Define
\begin{eqnarray}
  \label{eq:definitions}
  S &=& \left\{\pi\colon\textstyle{\sum_i\pi_i=1}\right\},\\
  S_i(q) &=& \left\{\pi\in{}S\colon{}i=\arg\max\{\pi_i/q_i\}\right\},\\
  A_i(q) &=& \int_{S_i(q)}\pi\Pr(D\mid{}\pi)\Pr(\pi)\,d\pi,\\
  C &=& \int_{S}\Pr(D\mid{}\pi)\Pr(\pi)\,d\pi.\label{eq:C}
\end{eqnarray}
The strategy given by problem~(\ref{eq:fopt}) in equivalent to
\begin{equation}
  \label{eq:fopt2}
  \min \sum_iq_i \quad\mbox{subject to}\quad 
  \sum_i\frac{A_i(q)}{q_i} = C.
\end{equation}
Observe that even if $\Pr(D\mid{}\pi)\Pr(\pi)$ has a complex form, if it can be
computed by Monte Carlo methods, then all the $A_i(q)$ and $C$ can be
computed simultaneously. The main difficulty is that the constraints have
to recomputed for each~$q$.

\subsection{Forced bets of entire wealth; logarithmic utility}\label{sec:allwealth}
We now consider the situation where the game rules require the client to
bet their entire wealth. This is a natural mathematical extension of
infinitesimal forced bets. With probabilistic odds this game leads to the
Kelly betting strategy of distributing the client's entire wealth as bets
proportional to the probabilities of events. This strategy maximises the
growth rate of wealth, or equivalently, a logarithmic utility function of
wealth. However, as \citet{Kelly:betting} explains, when the odds are not
probabilistic, then the client should not bet their entire wealth, only a
fraction of it. We will insist that the client is forced to bet their
entire wealth. The unfairness of this restriction is balanced by the
information advantage of the client.  If the forecaster tries to make the
odds as probabilistic as possible, then the odds can still interest the
client.

Once again, suppose there is a complete set of mutually exclusive
events~$E_i$, $i=1,\dots,m$ and a $m\times{}m$ game matrix $\cal{}G$ with
$G_{ij}=(\delta_{ij}/q_i)-1$. The client plays nature and is required to
bet their entire wealth at each play. Suppose that after $n$ plays of this
game, the client has total wealth~$W_n$, and at each play the client
distributes their current wealth as bets in the proportion $p_i$ on
event~$E_i$. Assume nature acts as though the events are chosen randomly,
with $E_i$ chosen with probability~$\pi_i$. Then the expected wealth of
the client after $n$ plays, given an initial wealth~$W_0$, is
\begin{displaymath}
  \label{eq:wealth}
  E(W_n) = \left(\prod_{i=1}^m\left(\frac{p_i}{q_i}\right)^{\pi_i}\right)^nW_0
\end{displaymath}
Consequently, the average rate of growth of wealth is
\begin{displaymath}
  \label{eq:growth}
  \lim_{n\to\infty}\frac{1}{n}\log(W_n) = \sum_{i=1}^m\pi_i\log(p_i/q_i).
\end{displaymath}
It is easily shown that the maximum growth rate occurs when $p_i=\pi_i$.
This is the optimal Kelly betting strategy for a client that knows
nature's~$\pi_i$. Hence, define $G(\pi,q)=\sum_{i=1}^m\pi_i\log(\pi_i/q_i)$.

An appropriate forecaster strategy is to set the closest to probabilistic
odds so that expected rate of growth of wealth of the client is zero, that
is,
\begin{equation}
  \label{eq:oddsgrowth}
  \min \sum_{i=1}^mq_i \quad\mbox{subject to}\quad
  \int_SG(\pi,q)\Pr(D\mid{}\pi)\Pr(\pi)\,d\pi=0.
\end{equation}
Equation~(\ref{eq:oddsgrowth}) implies that the logarithm of the client's
wealth is a martingale, that is, $E(\log{}W_n\mid{}D,q) =
\log{}W_{n-1}+E(G\mid{}D,q) = \log{}W_{n-1}$.

Remarkably the optimisation of equation~(\ref{eq:oddsgrowth}) has an
explicit form for the solution, which means that in some instances closed
form solutions are possible.  Define, using $C$ as in
equation~(\ref{eq:C}),
\begin{eqnarray}
  \label{eq:barHpi}
  \bar{H} &=& \frac{1}{C}\int_S\sum_{i=1}^m\pi_i\log(\pi_i)\Pr(D\mid{}\pi)\Pr(\pi)\,d\pi,\\
  \bar\pi_i &=& \frac{1}{C}\int_S\pi_i\Pr(D\mid{}\pi)\Pr(\pi)\,d\pi,\\
  {\alpha} &=& 
  \bar{H}-\sum_{i=1}^m\bar\pi_i\log(\bar\pi_i).
\end{eqnarray}
The constraint $E(G\mid{}D,q)=0$ in (\ref{eq:oddsgrowth}) is equivalent to
$\sum_{i=1}^m\bar\pi_i\log(q_i)=\bar{H}$. By straight forward application
of Lagrange multipliers it can be shown that $q_i = \bar\pi_i e^{\alpha}$
solves the required optimisation~(\ref{eq:oddsgrowth}).  The explicit form
of the solution, in terms of constants that are obtained from integrals,
means that in some instances closed form solutions are possible. A valid
interpretation of this optimisation result is that the inflation exponent
$\alpha$ is the discrepancy between the expected entropy (or
information)~$-\bar{H}$ and the entropy implied by the expected
probabilities~$\bar\pi$.

%%%%%%%%%%%%%%%%%%%%%%%%%%%%%%%%%%%%%%%%%%%%%%%%%%%%%%%%%%%%%%%%%%%%%%%%%%%

\section{Examples of odds forecasting}\label{sec:models}
To illustrate the computation of odds we consider two situations
forecasting binary events. The first situation uses a \emph{frequency
  model}, which requires only knowledge of the frequency of past events.
The aim of the forecaster is to provide odds on a future event. (This is
the urn game.) The second situation uses a \emph{Gaussian model}, where
the available data is a finite collection of scalar measurements assumed
to be drawn from a Gaussian distribution.  The aim of the forecaster is to
provide odds on a future measurement being below some threshold. We will
actually consider four situations because we will compute the odds for
both linear and logarithmic utility.

\subsection{Frequency model with linear utility}\label{sec:freq}
Consider the situation where there are two events $E$ and its
complement~$E'$, where Nature selects $E$ with probability~$\pi$ and $E'$
with probability~$\pi'=1-\pi$. A frequency model assumes that all
realizations of the events are independent, so only the frequency of the
event provides any information about~$\pi$. The forecaster is required to
assign odds $q$ and~$q'$ to events $E$ and~$E'$.

Suppose the event~$E$ has been observed to occur~$x$ times in~$n$
realizations.  Under the frequency model
$\Pr(\pi\mid{}x,n)\propto\Pr(x\mid{}\pi,n)\Pr(\pi)=C^n_x\pi^x(1-\pi)^{n-x}\Pr(\pi)$.
Suppose the forecaster has no prior information on~$\pi$, and so assumes
a uniform prior~$\Pr(\pi)=1$, that is, assumes all values of $\pi$ are
equally likely.

Following equation~(\ref{eq:fopt2}) the constraint on $q$ and~$q'$ to
obtain $E(V\mid{}x,n,q,q')=0$ can be expressed in terms of beta and incomplete
beta functions,
\begin{equation}
  \label{eq:qqfreq}
  \frac{1}{q'}\beta_\frac{q}{q+q'}(x+1,n-x+2) +
    \frac{1}{q}\beta_\frac{q'}{q+q'}(n-x+1,x+2) = \beta(x+1,n-x+1),
\end{equation}
where $\beta_x(a,b)=\int_0^xt^{a-1}(1-t)^{b-1}\,dt$ and
$\beta(a,b)=\beta_1(a,b)$.  Defining $s=q+q'$ and $q=ps$, so $q'=(1-p)s$,
then solving equation~(\ref{eq:qqfreq}) for~$s$, will conveniently
transform problem~(\ref{eq:fopt2}) into a one-dimensional problem,
\begin{equation}
  \label{eq:fopt3}
  \min_{0\leq{p}\leq1} s = \frac{\frac{1}{1-p}\beta_p(x+1,n-x+2) +
    \frac{1}{p}\beta_{1-p}(n-x+1,x+2)}{\beta(x+1,n-x+1)}.
\end{equation}
This problem is easily solved numerically using Brent's method or
similar~\citep{Press-etal:numerical-recipes}.

\begin{figure}[htbp]
  \centering
  \begin{tabular}{ll}
    \begin{tabular}{|r|c|c|c|}
      \hline
      $x/n$ & $q$ & $q+q'$\\
      \hline
      $0/1$ &    0.556 &    1.411 \\
      $1/1$ &    0.855 &    1.411 \\
      \hline
      $0/2$ &    0.448 &    1.347 \\
      $1/2$ &    0.687 &    1.375 \\
      $2/2$ &    0.900 &    1.347 \\
      \hline
      $0/3$ &    0.378 &    1.302 \\
      $1/3$ &    0.578 &    1.336 \\
      $2/3$ &    0.758 &    1.336 \\
      $3/3$ &    0.924 &    1.302 \\
      \hline
      $0/4$ &    0.329 &    1.268 \\
      $1/4$ &    0.501 &    1.303 \\
      $2/4$ &    0.656 &    1.313 \\
      $3/4$ &    0.803 &    1.303 \\
      $4/4$ &    0.939 &    1.268 \\
      \hline
    \end{tabular}
    &
    \raisebox{-4cm}{\includegraphics[width=11.5cm]{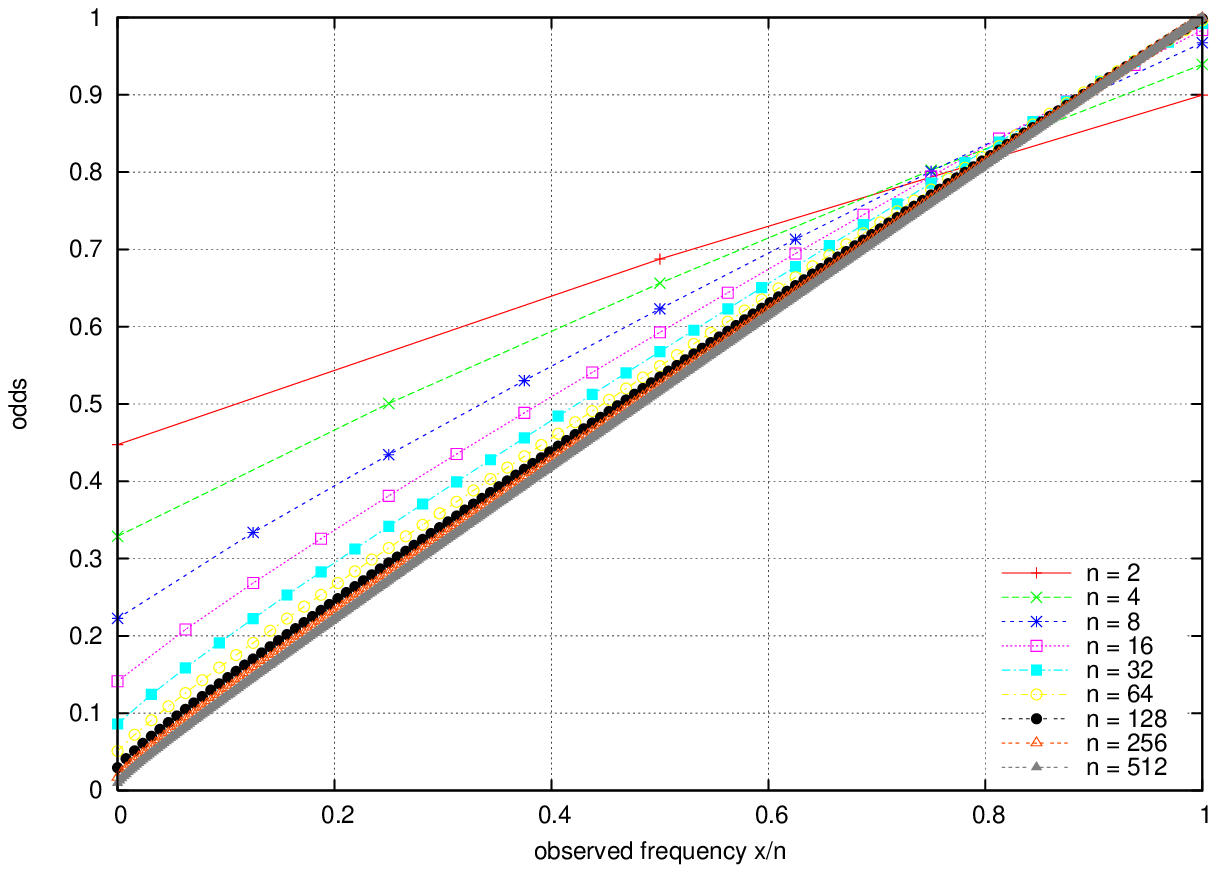}}
  \end{tabular}
  \caption{%{}
    Odds for frequency model with linear utility function. The table shows
    for frequencies $x/n$ of event~$E$, and the values of $q$ and $q+q'$,
    as computed from~(\ref{eq:fopt3}). The value of $q'$ can be obtained,
    by symmetry, from the frequency~$(n-x)/n$.  Values for larger~$n$ are
    plotted in fig.~\ref{fig:odds}. Only marked points are realizable
    values of the frequency.}
  \label{fig:odds}
\end{figure}

\begin{figure}[htbp]
  \centering
  \includegraphics[width=12cm]{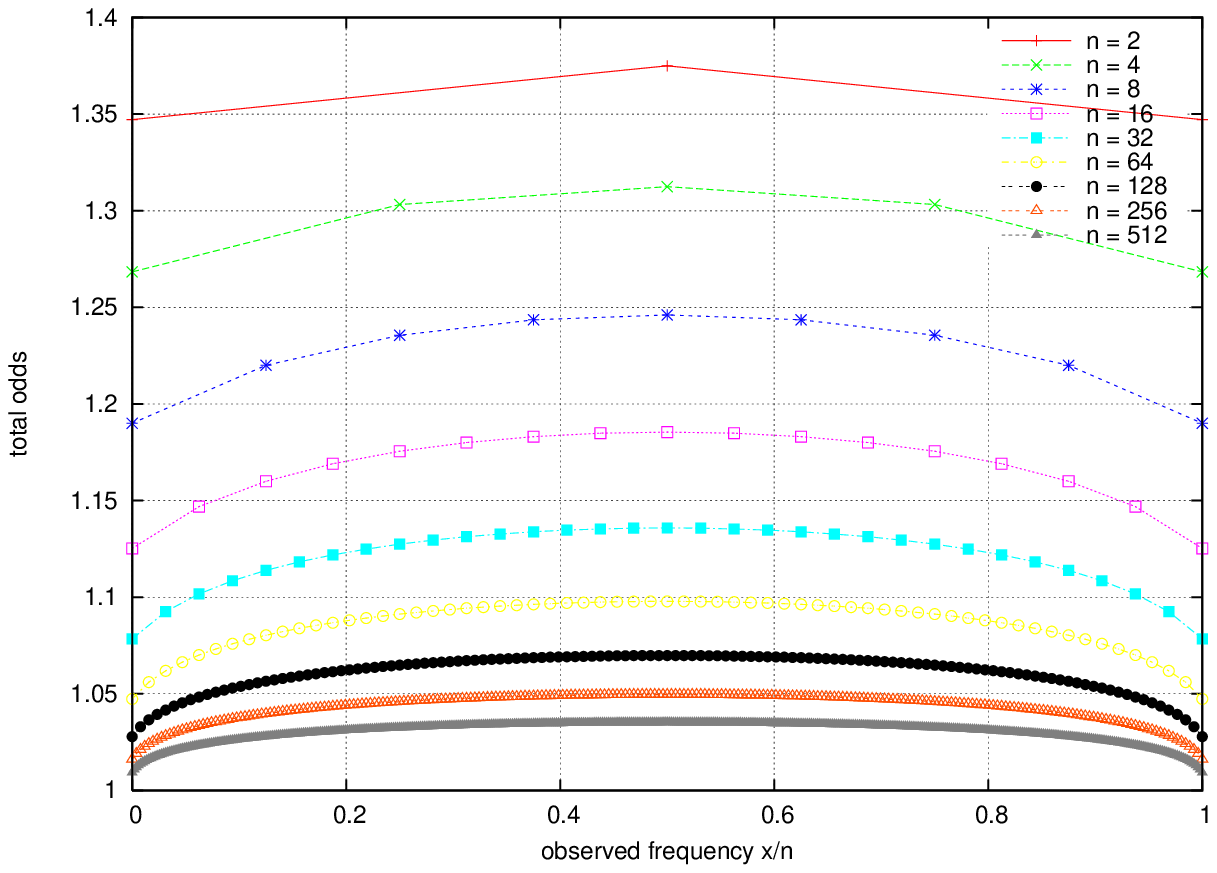}
  \caption{%{}
    The total odds $s=q+q'$ for frequency model with linear utility
    function. Only marked points are realizable values of the frequency.}
  \label{fig:totalodds}
\end{figure}

Figure~\ref{fig:odds} shows computed odds for various observed frequencies
for a small number of observations in a table, and for progressively
larger numbers of observations in the graph.  Figure~\ref{fig:totalodds}
shows the total odds $s=q+q'$. A number of interesting, but not
unexpected, facts can be seen. For a small number of observations the odds
are far from probabilistic, but they become more probabilistic as the
number of observations increase. Furthermore, the odds deviate most from a
probability for the event with a low frequency count, which is consistent
with the base-rate effect. Observe that meaningful odds are given when
$x=0$ or~$n$, and that $q\to{}x/n$ and $q'\to{}(n-x)/n$ as $n\to\infty$.
Furthermore, we find that for a fixed ratio $x/n$,
$q+q'\to{}1+\frac{c}{\sqrt{n}}$ for some constant~$c$ as $n\to\infty$.

%%%%%%%%% figures

\subsection{Frequency model with logarithmic utility}
This situation allows an essentially closed form solution for the odds.
Using equation~(\ref{eq:barHpi}) with
$\Pr(x,n\mid{}\pi)=C^n_x\pi^x(1-\pi)^{n-x}$ and $\Pr(\pi)=1$, obtains
\begin{eqnarray*}
  q &=& \left(\frac{x+1}{n-x+1}\right)^{\frac{n-x+1}{n+2}} e^\psi,\qquad
  q' = \left(\frac{n-x+1}{x+1}\right)^{\frac{x+1}{n+2}} e^\psi,\\
  \psi &=& 
  \left(\frac{x+1}{n+2}\right)H_{x+1}+\left(\frac{n-x+1}{n+2}\right)H_{n-x+1}-H_{n+2}
  %\frac{1}{n+2}\left((n-x+1)\psi(n-x+2)+(x+1)\psi(x+2)\right)-\psi(n+3)
\end{eqnarray*}
where we use the harmonic numbers $H_k =\sum_{i=1}^k\frac{1}{i}$.

\begin{figure}[htbp]
  \centering
  \begin{tabular}{ll}
    \begin{tabular}{|r|c|c|c|}
      \hline
      $x/n$ & $q$ & $q+q'$\\
      \hline
      $0/1$ &    0.382 &    1.146 \\
      $1/1$ &    0.764 &    1.146 \\
      \hline
      $0/2$ &    0.277 &     1.11 \\
      $1/2$ &    0.558 &    1.116 \\
      $2/2$ &    0.832 &     1.11 \\
      \hline
      $0/3$ &    0.217 &    1.087 \\
      $1/3$ &    0.438 &    1.094 \\
      $2/3$ &    0.656 &    1.094 \\
      $3/3$ &     0.87 &    1.087 \\
      \hline
      $0/4$ &    0.179 &    1.073 \\
      $1/4$ &    0.359 &    1.078 \\
      $2/4$ &     0.54 &    1.079 \\
      $3/4$ &    0.719 &    1.078 \\
      $4/4$ &    0.894 &    1.073 \\
      \hline
    \end{tabular}
    &
    \raisebox{-4cm}{\includegraphics[width=11.5cm]{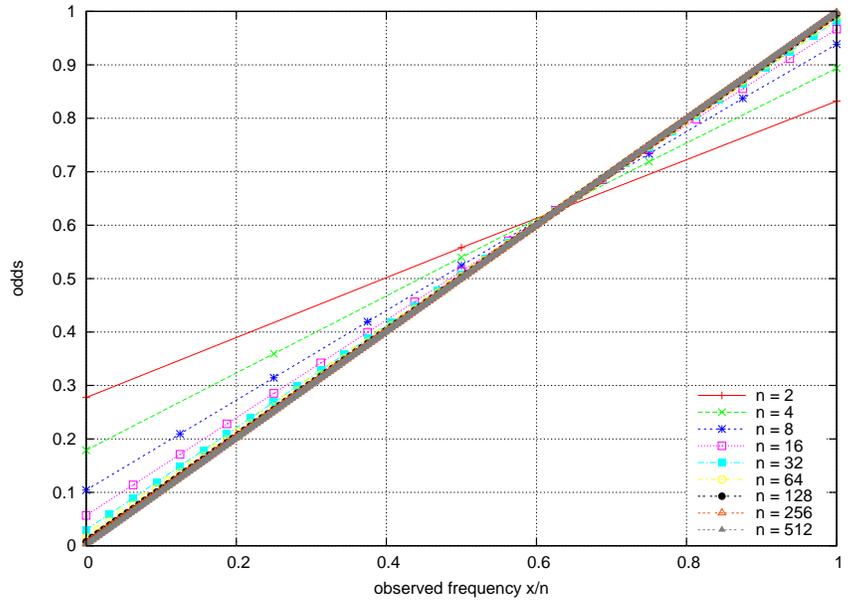}}
  \end{tabular}
  \caption{%{}
    Odds for frequency model with logarithmic utility function. The table
    shows for frequencies $x/n$ of event~$E$, and the values of $q$ and
    $q+q'$, as computed from the closed form formula. The value of $q'$ can be
    obtained, by symmetry, from the frequency~$(n-x)/n$.  Values for
    larger~$n$ are plotted in fig.~\ref{fig:odds}. Only marked points
    are realizable values of the frequency.  }
  \label{fig:oddsgg}
\end{figure}

\begin{figure}[htbp]
  \centering
  \includegraphics[width=12cm]{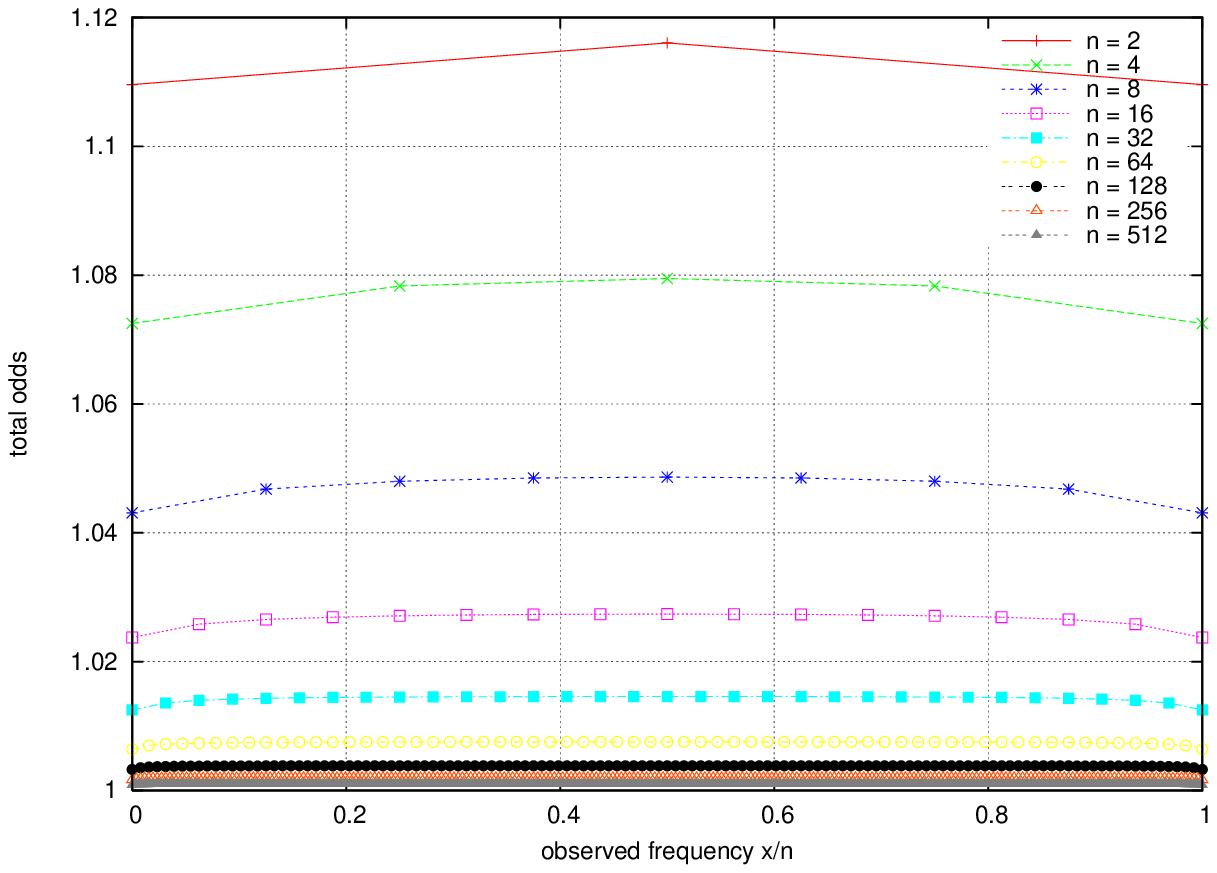}
  \caption{%{}
    The total odds $s=q+q'$ for frequency model with logarithmic utility
    function. Only marked points are realizable values of the frequency.}
  \label{fig:totaloddsgg}
\end{figure}

Figure~\ref{fig:oddsgg} provides a table of computed values of the odds
for small values of~$n$ and graphs the odds for larger values. These
computed odds should be compared with fig.~\ref{fig:odds}.
Figure~\ref{fig:totaloddsgg} shows how the total odds~$q+q'$ varies with
number of observations, which should be compared with
fig.~\ref{fig:totalodds}. It is observed that the total odds have a
smaller excess over unity, and vary less with the observed fraction. This
is not surprising, because a client who aims to maximise their rate of
growth of wealth is less likely to exploit forecast errors of low
probability events. Observe also that the odds converge much faster in
this situation, at a rate of $1/n$, as opposed to $1/\sqrt{n}$ in the
infinitesimal bets case. (The author also has closed form expressions for
the odds in this case of an arbitrary number of events, which will
discussed elsewhere.)

\subsection{Gaussian model with linear utility}\label{sec:gaussian}
Now consider a situation where the forecaster is given scalar observations
$D=\{x_1,\dots,x_n\}$, with statistics
\begin{eqnarray*}
  \label{eq:statistics}
  \hat\mu&=&\frac{1}{n}\sum_ix_i\\
  \hat\sigma^2&=&\frac{1}{n}\sum_i(x_i-\hat\mu)^2
  =\frac{1}{n}\sum_ix_i^2-\hat\mu^2.
\end{eqnarray*}
Furthermore, the forecaster has reason to model these as observations of a
random variable~$X$ with a Gaussian distribution~$N(\mu,\sigma^2)$, for
some unknown $\mu$ and~$\sigma$; although the forecaster may have
additional prior belief~$\Pr(\mu,\sigma)$ in the values of $\mu$
and~$\sigma$. The client requires an odds forecast for the events
$E=\{X\le{}x\}$ and $E'=\{X>x\}$ for some fixed~$x$.

Since we are assuming a perfect model class (to avoid third order
uncertainty issues), then the true probability~$\pi$ of the event~$E$ is
\begin{displaymath}
  \label{eq:normal}
  \pi=\Pr(X\le{}x)=
  \int_{-\infty}^x\frac{1}{\sqrt{2\pi\sigma^2}}
  e^{-(t-\mu)^2/2\sigma^2}\,dt
  = \Phi\left(\frac{x-\mu}{\sigma}\right).
\end{displaymath}
Whereas the optimisation problem~(\ref{eq:fopt}) that obtains the odds is
formulated in terms of integrals over the probabilities of the events, in
this situation we see these probabilities are determined entirely from
$x$, $\mu$ and $\sigma$. Consequently, it is more appropriate to
reformulate the optimisation~(\ref{eq:fopt}) it terms of integrals over
$\mu$ and $\sigma$. This requires reformulating equations
(\ref{eq:V}) and (\ref{eq:EV}).

When the forecaster assigns odds $q$ and~$q'$ to $E$ and~$E'$
respectively, then in the worst case where the client knows $\mu$
and~$\sigma$ (and hence $\pi$) the optimal strategy of the client is to
bet on
\begin{eqnarray*}
  \label{eq:nclient}
  E &\mbox{if}& \frac{q}{q+q'}<\Phi\left(\frac{x-\mu}{\sigma}\right),\\
  E'&\mbox{if}& \frac{q'}{q+q'}<1-\Phi\left(\frac{x-\mu}{\sigma}\right).
\end{eqnarray*}
and the average pay-out to client is then
\begin{displaymath}
  \label{eq:npayout}
  V(q,q'\mid{}x,\mu,\sigma) = \max
  \left\{
    \frac{1}{q}\Phi\left(\frac{x-\mu}{\sigma}\right),
    \frac{1}{q'}\left(1-\Phi\left(\frac{x-\mu}{\sigma}\right)\right)
  \right\}-1.
\end{displaymath}

Given their model the forecaster will assert that
$\Pr(\pi\mid{}D)\propto\Pr(D\mid{}\mu,\sigma)\Pr(\mu,\sigma)$ where
\begin{eqnarray*}
  \label{eq:pm}
  \Pr(D\mid{}\mu,\sigma) &=& 
  \prod_{i=1}^n\frac{1}{\sqrt{2\pi\sigma^2}}e^{-(x_i-\mu)^2/2\sigma^2}\\
  &=& \left(\frac{1}{\sqrt{2\pi\sigma^2}}\right)^n
  e^{-n(\hat\sigma^2+(\hat\mu-\mu)^2)/2\sigma^2}.
\end{eqnarray*}

Just as in section~\ref{sec:freq} we  once again have a situation
involving binary events where there is advantage in defining $s=q+q'$,
$q=ps$, $q'=(1-p)s$, so that problem~(\ref{eq:fopt}) can be transformed to
\begin{displaymath}
  \label{eq:fopt4}
  \min_{0\leq{p}\leq1}s = 
  \frac{1}{C}\left(\frac{A(p,x)}{1-p}+\frac{B(p,x)}{p}\right),
\end{displaymath}
where
\begin{eqnarray*}
  \label{eq:cc}
  A(p,x) &=& \int_0^\infty\int_{M(p,x,\sigma)}^\infty
  \Phi\left(\frac{x-\mu}{\sigma}\right)
  F(\mu,\sigma)
  \,d\mu\,d\sigma,\\
  B(p,x) &=& \int_0^\infty\int_{-\infty}^{M(p,x,\sigma)}
  \left(1-\Phi\left(\frac{x-\mu}{\sigma}\right)\right)
  F(\mu,\sigma)
  \,d\mu\,d\sigma,\\
  C &=& \int_0^\infty\int_{-\infty}^{\infty}
  F(\mu,\sigma)
  \,d\mu\,d\sigma,\\
  M(p,x,\sigma) &=& x-\sigma\Phi^{-1}(p),\\
  F(\mu,\sigma) &=& \Pr(D\mid{}\mu,\sigma)\Pr(\mu,\sigma).
\end{eqnarray*}

The odds can be calculated once the prior $\Pr(\mu,\sigma)$ has been
assigned, however, the assignment of the prior requires a little care. It
is not unreasonable to assume no knowledge of $\mu$, and hence place a
uniform prior on~$\mu$. On the other hand, one cannot assume a uniform
prior for~$\sigma$, because the integrals will diverge. Taking a uniform
prior for $\sigma$ on $(0,R)$, then one finds that $s\to1$ and $p\to1/2$
as $R\to\infty$. Essentially a model that freely allows arbitrarily large
variances cannot provide useful forecasts. Consequently, one is required
to specify a prior for $\sigma$ with a tail that thins sufficiently
rapidly.

\begin{figure}[htbp]
  \centering
  \includegraphics[width=12cm]{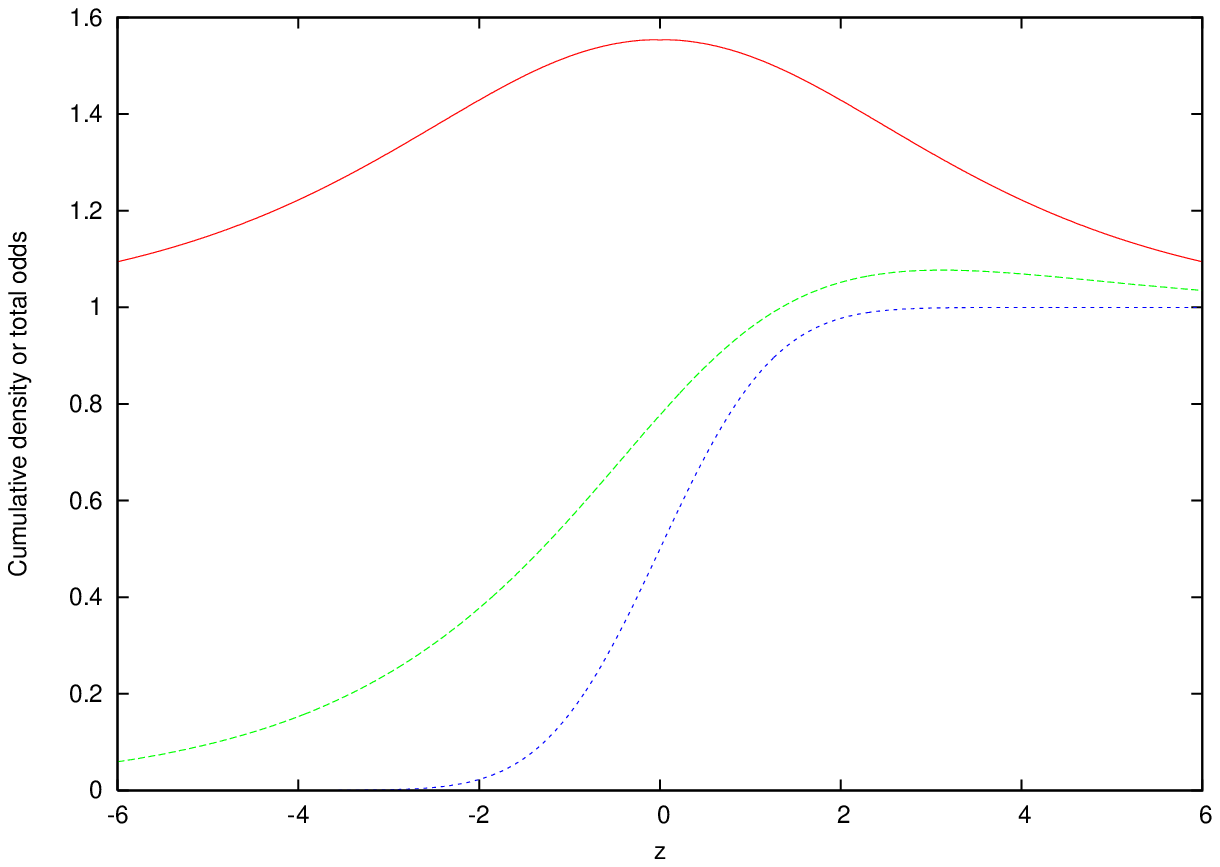}
  \caption{%{}
    Odds~$q$ (dashed) and total odds $s=q+q'$ (solid) for Gaussian model
    with linear utility function, here computed for the generic situation
    $\hat{\mu}=0$ and $\hat\sigma=1$, prior
    $\Pr(\mu,\sigma)=\sigma{}e^{-\sigma^2}$, and events $E=\{Z\le{}z\}$
    and $E'=\{Z>z\}$.  The standard normal cumulative distribution
    function is also plotted for comparison (dotted).}
  \label{fig:godds}
\end{figure}

Figure~\ref{fig:godds} shows numerically computed odds for a situation
where $\hat{\mu}=0$ and $\hat{\sigma}=1$. The prior used was
$\Pr(\mu,\sigma) = \sigma e^{-\sigma^2/2}$, which is a $\chi$-distribution
of two degrees of freedom. This prior was used because we want
$\hat{\sigma}=1$ to be a typical observed value. The computed odds were
similar for half-normal, F-distributions and other similar distributions
where $\hat{\sigma}=1$ is a typical observed value.
By a shift and scale the computed odds shown in fig.~\ref{fig:godds}
provide a general solution for arbitrary $\hat{\mu}$ and $\hat{\sigma}$.
Let $x_i=\hat{\mu}+\hat{\sigma}z_i$ so that the $z_i$ have mean zero and
variance one. These $z_i$ should be modelled as observations of a random
variable $Z=(X-\hat{\mu})/\hat{\sigma}$.  The events of interested are now
expressed as $E=\{Z\le{}z\}$ and $E'=\{Z>z\}$, where
$x=\hat{\mu}+\hat{\sigma}z$. Hence, if $q(z)$ and $s(z)$ represent the
odds and total odds for the generic case shown in fig.~\ref{fig:godds},
then the odds in the general case are obtained using
$z=(x-\hat{\mu})/\hat{\sigma}$.

Figure~\ref{fig:godds} shows that the odds on~$E$ exceed one for $z$
larger than about 1. Our interpretation of this is that given the
information available to the forecaster, the event~$E$ is so likely that
offering odds favourable to the client would be loss making to the
forecaster. The client is therefore offered odds so that to bet on~$E$
makes a consistent small loss with an rare loss of $-1$, and to bet
on~$E'$ makes a consistent loss of $-1$ and a rare win.

\subsection{Gaussian model with logarithmic utility}
Computation of odds in this situation follow in a similar fashion to the
previous section, in that the integrals~(\ref{eq:barHpi}) are reformulated
in terms of integrals over $\mu$ and $\sigma$. Thus, using the same
notation as the previous sections, the odds are $q=\bar\pi e^\alpha$ and
$q'=(1-\bar\pi) e^\alpha$ where
\begin{eqnarray*}
  \label{eq:barHpi2}
  \bar{H} &=& 
   \frac{1}{C}\int_0^\infty\int_{-\infty}^\infty
  \Phi((z-\mu)/\sigma)\log(\Phi((z-\mu)/\sigma))F(\mu,\sigma)\,d\mu\,d\sigma,\\
  && +
   \frac{1}{C}\int_0^\infty\int_{-\infty}^\infty
  (1-\Phi((z-\mu)/\sigma))\log(1-\Phi((z-\mu)/\sigma))F(\mu,\sigma)\,d\mu\,d\sigma,\\
  \bar\pi &=&  \frac{1}{C}\int_0^\infty\int_{-\infty}^\infty
  \Phi((z-\mu)/\sigma)F(\mu,\sigma)\,d\mu\,d\sigma,\\
  {\alpha} &=& 
 \bar{H} -\bar\pi\log(\bar\pi)-(1-\bar\pi)\log(1-\bar\pi)
\end{eqnarray*}

\begin{figure}[htbp]
  \centering
  \includegraphics[width=12cm]{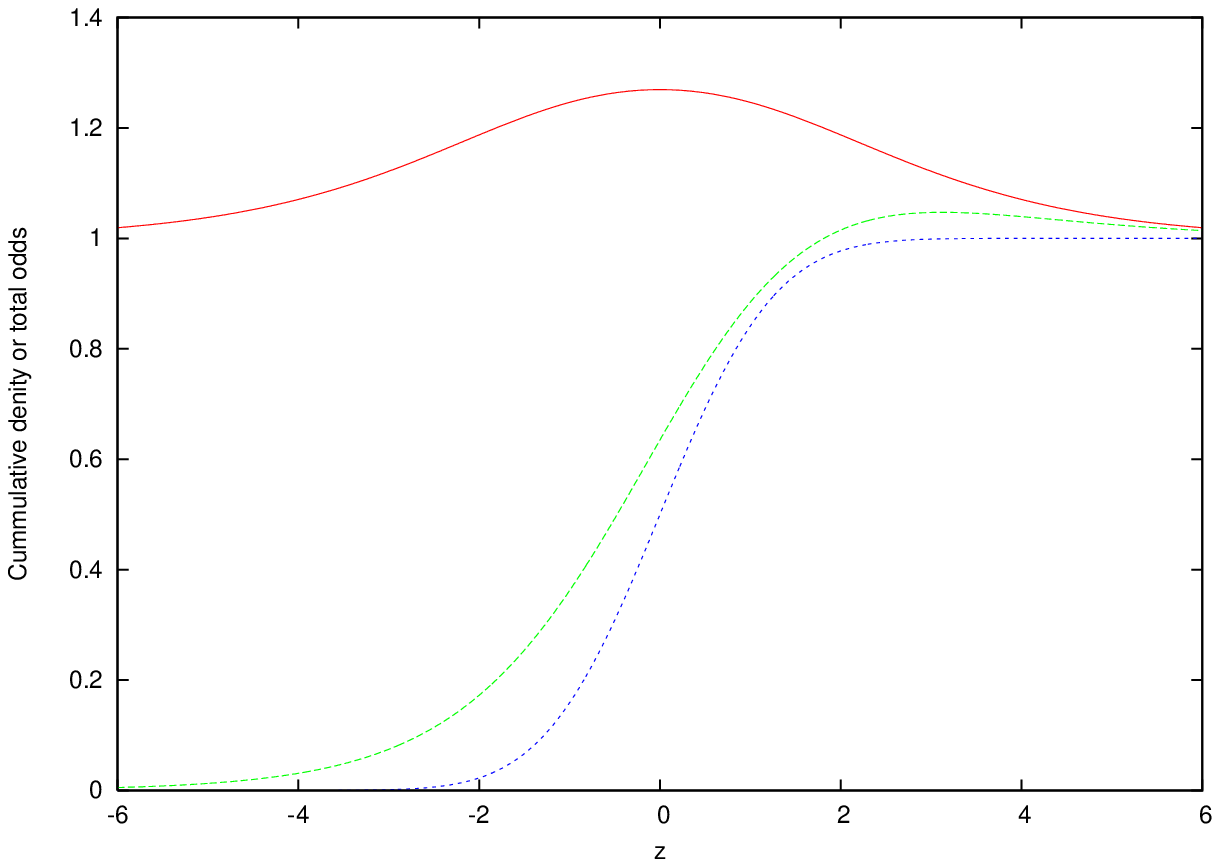}
  \caption{%{}
    Odds~$q$ (dashed) and total odds $s=e^\alpha$ (solid) for Gaussian
    model with logarithmic utility function, here computed for the generic
    situation $\hat{\mu}=0$ and $\hat\sigma=1$, prior
    $\Pr(\mu,\sigma)=\sigma{}e^{-\sigma^2}$, and events $E=\{Z\le{}z\}$
    and $E'=\{Z>z\}$.  The standard normal cumulative distribution
    function is also plotted for comparison (dotted).}
  \label{fig:gmodds}
\end{figure}

Figure~\ref{fig:gmodds} shows the computed odds for a situation where
$\hat{\mu}=0$ and $\hat{\sigma}=1$, and the prior $\Pr(\mu,\sigma) =
\sigma e^{-\sigma^2/2}$. These odds should be compared with the linear
utility situation shown in fig.~\ref{fig:godds}. In comparison it is
seen that the logarithmic utility function gives odds with smaller
excess~$s$ and less weight attached to the large~$z$ values. this is
similar to the frequency model, for the same reasons.

%%%%%%%%%%%%%%%%%%%%%%%%%%%%%%%%%%%%%%%%%%%%%%%%%%%%%%%%%%%%%%%%%%%%%%%%%%%
\section{Investment and loss mitigation}\label{sec:investment}

The client who has featured thus far in our analysis is more accurately
described as a \emph{speculative client}, whose primary goal is to profit
from inadequacies of the forecaster's predictions. We now introduce the
\emph{invested client}, whose wealth is invested in some venture whose
profit, costs, and losses are determined in part by the outcomes of the
events. Think here of the road-gritter, or horticulturist who is concerned
with the possibility of freezing temperatures. The invested client has no
desire or ability to challenge the forecaster's skill, rather they wish to
use the forecasts to mitigate their losses. The invested client can do
this by betting against the forecaster; a bet is essentially an insurance
policy. In this section we analyse how an invested client should bet and
show that such bets are beneficial to the invested client despite the
forecaster's odds being non-probabilistic.  Furthermore, we will see that
the client deals only with the forecaster's odds, they do not try to
normalise the odds to obtain probability ``estimates''; the odds contain
all the information the invested client needs.

\subsection{Simple investment}
Consider a situation where the return on a client's investment is
influenced by whether an event~$E$, or its complement~$E'$, occurs. Let
$R$ be the return on the investment when the event~$E$ occurs and $R'$ the
return on $E'$.  Let $\pi$ and~$\pi'$, where $\pi+\pi'=1$, be the
probabilities of the events $E$ and~$E'$ respectively, under assumption
that Nature selects the events at random. The expected return on the
investment is $R\pi+R'\pi'$.

Now suppose the forecaster provides odds $q$ and~$q'$ for the events~$E$
and~$E'$. The game matrix for betting of these events
is~(\ref{eq:game2x2}), in terms of the pay-outs $P=(1/q)-1$ and
$P'=(1/q')-1$.  If the client places bets with the forecaster of
$\lambda\ge0$ on the event~$E$ and $\lambda'\ge0$ on the complement
event~$E'$, then the expected return to the client is
\begin{displaymath}
  \label{eq:return}
  (R+\lambda{}P-\lambda')\pi + (R'-\lambda+\lambda'P')\pi'.
\end{displaymath}
If the client chooses the bets $\lambda$ and $\lambda'$ so that
$R+\lambda{}P-\lambda' = R'-\lambda+\lambda'P'$, then the client's return
is fixed and independent of $\pi$ and~$\pi'$, indeed independent of which
event occurs. Under this condition
\begin{displaymath}
  \label{eq:equal2}
  R-R' = (P'+1)\lambda' - (P+1)\lambda = 
  \frac{\lambda'}{q'}-\frac{\lambda}{q},
\end{displaymath}
and so, if $R>R'$, then the client should bet $\lambda=0$ and
$\lambda'=q'(R-R')$ for a return of $(1-q')R+q'R'$, and if $R<R'$, then
the client should bet $\lambda=q(R'-R)$ and $\lambda'=0$ for a return of
$qR+(1-q)R'$. 

Since $q+q'>1$ the returns with bets are less than the expected return
without bets, the advantage to the client is that by placing bets the
return is guaranteed and the risk is transferred to the forecaster. In a
market of forecasters those that offer the odds closer to being
probabilities will attract more clients, but if their odds have not been
calculated along the lines we have described, then they will be almost
surely bankrupted.

\subsection{Mitigating losses}

Consider a situation where a client incurs a loss~$L$ if an event~$E$
occurs, but they can take an action at cost~$C$ which if taken results in
inclusive mitigated losses~$M$. This situation is represented by the
following game matrix.
\begin{displaymath}
  \label{eq:newgame}
  \begin{tabular}{cc}
    \raisebox{-0.5cm}{client} & 
    \begin{tabular}{c|cc|}
      \multicolumn{1}{c}{} & \multicolumn{2}{c}{nature}\\
      & $E$ & $E'$ \\
      \hline
      no action & $-L$ & $0$ \\
      action & $-M$ & $-C$ \\
      \hline
    \end{tabular}
  \end{tabular}  
\end{displaymath}
In this situation it is usual that $0<C<M<L$, although it can happen that
mitigating the losses includes a reward that more than covers costs of the
action so that $M<C$. In either case the optimal strategy for the client
is to take the action when $\pi>C/(L+C-M)$.

As before a forecaster provides odds on the events $E$ and~$E'$. Suppose
the client places bets $\lambda$ and $\lambda'$ on $E$ and~$E'$
respectively when no action is taken, and places bets $\mu$ and $\mu'$ when
the action is taken. The game matrix is the following.

\begin{displaymath}
  \label{eq:newgame2}
  \begin{tabular}{c|cc|}
    & $E$ & $E'$ \\
    \hline
    no action & $-L+\lambda{}P-\lambda'$ & $\lambda'{}P'-\lambda$ \\
    action & $-M+\mu{}P-\mu'$ & $-C+\mu'{}P'-\mu$ \\
    \hline
  \end{tabular}
\end{displaymath}

Following the analysis of the previous subsection there are three
possibilities:
\begin{enumerate}
\item[(a)] Take no action and place bets $\lambda=qL$ and $\lambda'=0$, which
  results in a fixed loss $-qL$.
\item[(b)] Given $C<M$, take the action and place bets $\mu=q(M-C)$ and
  $\mu'=0$, which results in a fixed loss $-qM-(1-q)C$.
\item[(c)] Given $M<C$, take the action and place bets $\mu=0$ and
  $\mu'=q'(C-M)$, which results in a fixed loss $-(1-q')M-q'C$.
\end{enumerate}

It follows that when $0<C<M<L$ the client always bets on the event
occurring and takes the action when $q>C/(L+C-M)$. In the situation where
$0<M<C<L$ the client bets on the event occurring when taking no action, and
bets on the event not occurring when taking the action, and takes the
action when $qL-(C-M)q'>M$. Once again if $q+q'>1$, then the losses with
bets are more than the expected losses without bets, but the losses are
fixed and the forecaster takes the risk.

The important point to note about these results, and those of the previous
subsection, is the client uses the odds $q$ and $q'$ directly, and does
not normalise these to obtain probability ``estimates'' $q/(q+q')$ and
$q'/(q+q')$ of the events. All the useful information to the client in
contained in the odds, and whether the clients refers to both $q$
and~$q'$, or just one of these values, depends on their circumstances. For
example, in the $0<C<M<L$ situation the decision to take the action is
based on $q$ alone, where as, in the $0<M<C<L$ situation it depends on
both $q$ and~$q'$, but not in a way that implies normalising them to
obtain probability ``estimates''.

%%%%%%%%%%%%%%%%%%%%%%%%%%%%%%%%%%%%%%%%%%%%%%%%%%%%%%%%%%%%%%%%%%%%%%%%%%%
\section{Ensemble Forecasting}\label{sec:ensembles}
Finally we consider an example of issuing of odds forecasts on temperature
variations using numerical ensemble weather predictions. A common event of
interest is whether the temperature at a locality will fall below
freezing.  These events are fairly rare, so we consider a related and more
general event of whether the temperature at a location at some set
lead-time falls more than a certain amount below the current temperature
at that station. These calculations are intended to provide an
illustration of odds forecasting using the results we have obtained so
far. We do not claim that the odds forecast we compute are the best or
most appropriate, indeed the results suggest otherwise. Odds based on
kernel density estimates, kernel dressing and the like could provide
better odds forecasts.

\subsection{Data preparation}
We use London Heathrow (Station 03772) as the locality and National Center
for Environmental Prediction GFS ensemble predictions~\citep{NCEP:GFS} for
constructing odds forecasts. The GFS model provides a temperature at the
station by interpolation of the temperatures of the global circulation
model's nearest grid-points to the location. Our first step is to prepare
suitable time series data including bias correction and temporal
interpolation.

The GFS 2-metre-temperature deterministic-forecast initialisation states
were used to calibrate GFS model temperatures with the station readings as
follows. Firstly, we compute a cubic spline of the time series of
deterministic-forecast initialisation states. Secondly, we compute using
the spline a model temperature time-series to pair with the station
temperature time-series. Finally, we fit by least squares the model
\begin{equation}
  \label{eq:transform}
  y=p_0(x)+p_1(x)\cos(2\pi{h}/24)+p_2(x)\sin(2\pi{h}/24),
\end{equation}
where $y$ is the station temperature, $x$ the splined model temperature,
$h$ the hour of the day, and each $p_j(x)$ is a cubic polynomial.  This is
a twelve parameter model. We used times series from the calendar year
2005. The residuals of this model had mean $6.03\times10^{-15}$, skewness
$0.0736$, kurtosis $0.9453$ and standard error of $1.1159$ degrees
centigrade.  Hence, this model provides a good transformation from model
state temperature to station temperature with an expected error that is
very nearly Gaussian and standard deviation of around $1.1159$ degrees
centigrade.

To obtain ensemble forecasts at any lead-time we compute a cubic spline to
the time series of each ensemble member, then apply our fitted
transform~(\ref{eq:transform}).

In the following these interpolated and transformed time-series will be
referred to as the \emph{adjusted control} and the \emph{adjusted ensemble
  forecasts}. Figure~\ref{fig:forecast} shows an example of the station
data, adjusted control and adjusted ensemble time-series.

\subsection{Odds forecasts}
The goal is to provide odds forecasts on whether the minimum daily
temperature at one to ten days lead-time will be less than 3 degrees
centigrade below the adjusted control temperature at the time of issuing
the forecasts.  In order to use the results we have already derived we
consider four forecasts.
The first forecast applies the odds from a frequency model, which is based
on counting number of adjusted ensemble members below the target
threshold.  The second forecast applies the odds from a Gaussian model,
which assumes the adjusted ensemble members have a Gaussian distribution.
The third forecast is a probabilistic forecast obtained by assuming the
adjusted ensemble members have a Gaussian distribution. The forth forecast
is intended to act as the ultimate-challenge client. It is not really a
forecast at all, but rather the predicted probability of station
temperature being below the threshold based on a Gaussian distribution
centred on the adjusted control with standard deviation equal to $1.1159$,
that is, error of the fitted residuals.

Forecasts one to three represent three competing forecasters, one and two
providing non-probabilistic odds, three providing probabilistic odds.
Forecast four is used by the ultimate-challenge client. Nature is always
taken to be the station temperature. Using the methods described in the
previous sections odds forecasts can be computed.
Figure~\ref{fig:oddsfcst} shows the odds computed using the time-series
data shown in fig.~\ref{fig:forecast} in the case of a linear utility.
Figure~\ref{fig:payout} shows the pay-out of the three forecasters when
playing against the ultimate-challenge client. Since the challenger is
using a linear utility they will bet on the temperature being below the
threshold if the forecaster's odds are less than the challenger's odds and
visa versa. Observe that the challenger loses more bets to forecaster~3,
who uses probabilistic odds, than to forecasters 1 and~2, who use
non-probabilistic odds, however, the winning bets on low probability
events against forecaster~3 are astronomical. On the other hand,
forecasters 1 and~2 appear to perform about the same.

We now consider the problem of forecasting whether the minimum temperature
at London Heathrow in the 24 hour period from 18:00UTC is more than 3
degrees centigrade below the control temperature when the forecasts are
issued. We will test the performance of the odds forecasts by computing
the total pay-out to the ultimate-challenge client for the first half of
2005 (183 days) for lead times up to ten days. Figures \ref{fig:forecast}
and~\ref{fig:oddsfcst} clearly demonstrate that the probabilistic odds we
have considered are not competitive with the non-probabilistic odds.
Certainly, large pay-outs result from assigning far too small probability
to low probability events, but what about the performance at other times?
To investigate this question we consider a new non-probabilistic odds
forecaster ($3'$) that uses the probabilistic odds but sets minimum
odds of $0.1$, which effectively caps pay-outs at 10. These capped
probabilities are a crude form of odds. Odds capping could be applied by
our forecasters 1 and~2, although in the test considered it makes little
difference.
Table~\ref{tab:odds} shows the total pay-out to the three non-probabilistic
odds forecasters 1, 2, and $3'$, on bets taken from the ultimate-challenge
client. If forecasters 1 and~2 are allowed odds caps of~$0.1$ is these
tests, then the cap only applies on around ten occasions and reduces the
larger pay-outs at lead times 5, 7, and 10 days to the level of their
neighbouring lead times.

A number of observations can be made about the results displayed in
table~\ref{tab:odds}. Recall that the aim of the forecasters is to have an
expected total pay-out of zero. It is immediately clear that odds capping
alone is not sufficient to make forecaster~$3'$ competitive with
forecasters 1 and~2 beyond a lead time of one or two days.  Forecasters 1
and~2 have similar performance although forecaster~2 is more successful.
In all cases it appears that performance decreases with increasing lead
time.

%%%%%%%%%%%%%%%%%%%%%%%%%%%%%%%%%%%%%%%%%%%%%%%%%%%%%%%%%%%%%%%%%%%%%%%%%%%
\section{Summary, discussion and conclusions}

We have confronted forecasters with the challenge that if they offer a
probability forecast, then they should be prepared to accept bets at the
odds these imply.  Unless a forecaster has a perfect forecast model, then
they would be unwise to accept this challenge, because any that do so will
almost surely be bankrupted by more informed bettors. We argue that
probability forecasts fail to account for higher order uncertainties, such
as model error. Our alternative is to offer non-probabilistic odds
forecasts obtain using an optimisation principle. The excess of the odds
reveals the forecaster's uncertainty about the model.  We have shown how
to compute odds forecasts in several situations, and illustrated how odds
forecasts could be used in investment, loss mitigation, and weather
forecasting.

\begin{figure}[htbp]
  \centering
  \includegraphics[width=12cm]{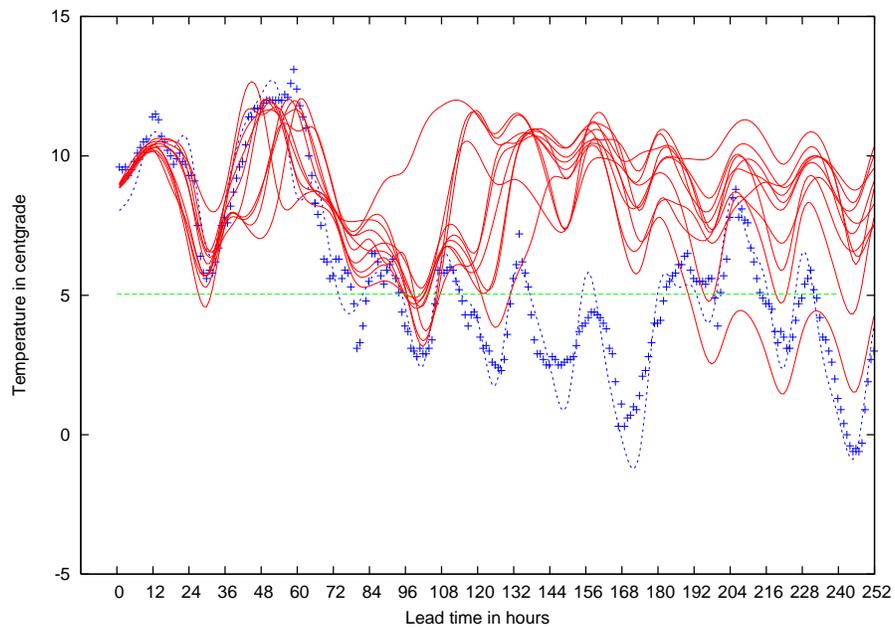}
  \caption{%{}
    Station temperature time series $(+)$, the adjusted control (dashed)
    and adjusted ensemble forecast time-series (solid). The horizontal
    line (dashed) shows the threshold line for which the odds forecasts
    are required. The forecasts shown were launched 00:00UTC day 40 of
    2005.}
  \label{fig:forecast}
\end{figure}

\begin{figure}[htbp]
  \centering
  \includegraphics[width=12cm]{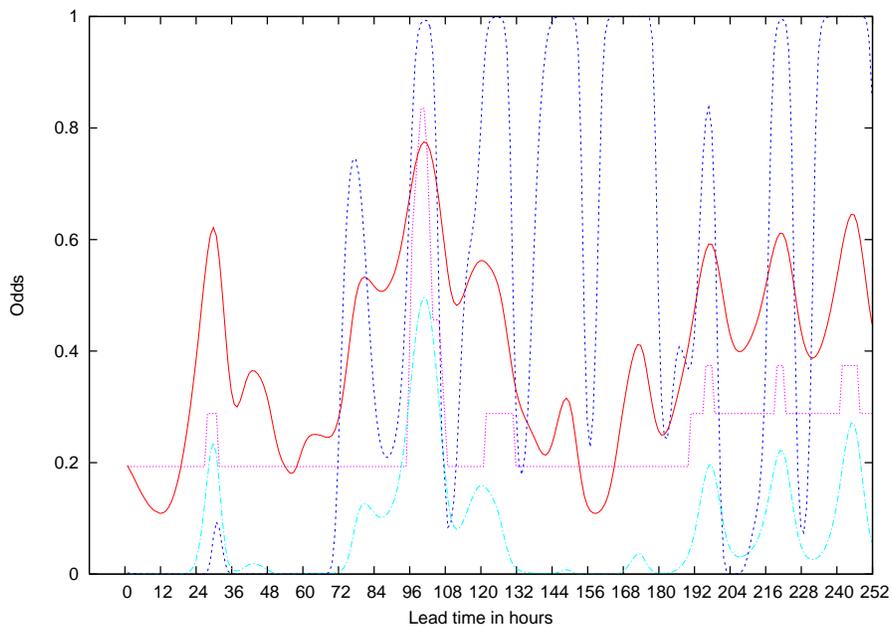}
  \caption{%{}
    Four computed odds forecasts (linear utility) for the data shown in
    figure~\ref{fig:forecast}. Frequency odds (dashed), Gaussian odds
    (solid), and probabilistic odds (dash-dotted).  Challenger (dashed)
    refers to the odds used by the ultimate-challenge client against whom
    the forecasters must compete.}
  \label{fig:oddsfcst}
\end{figure}

\begin{figure}[htbp]
  \centering
  \includegraphics[width=12cm]{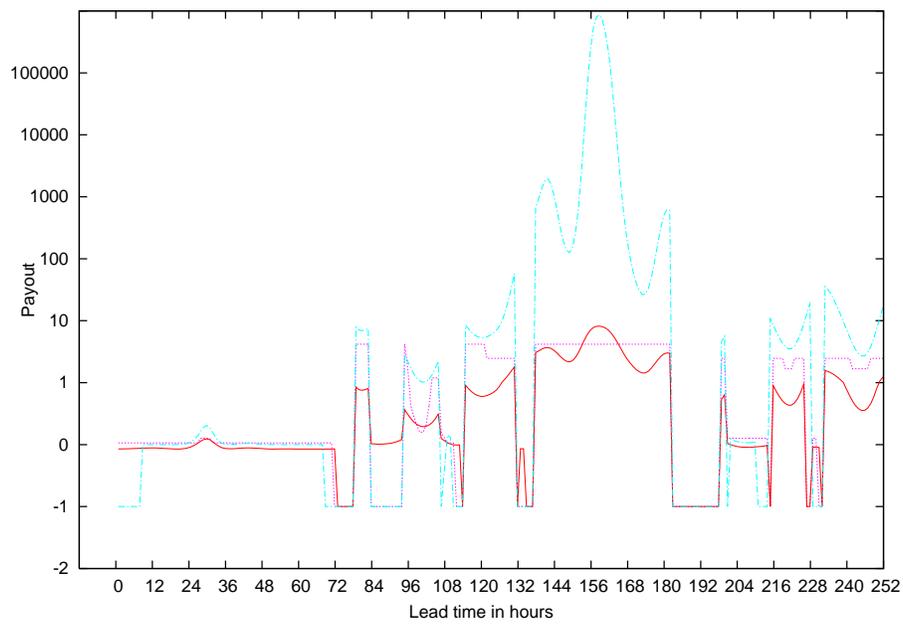}
  \caption{%{}
    The pay-outs of the three forecasters when faced with the
    ultimate-challenge client, for linear utility. Frequency odds
    (dashed), Gaussian odds (solid), and probabilistic odds (dash-dotted).
    Same data as shown is figures \ref{fig:forecast}
    and~\ref{fig:oddsfcst}.  The pay-outs are shown on a linear scale for
    pay-outs less than one, then transitions into a logarithm base 10
    scale above.}
  \label{fig:payout}
\end{figure}

There are gaps in our development and demonstrations, especially in the
application to ensemble weather forecasting. A gap occurs because we
have only shown how to compute odds under the assumption of a perfect
model class; this assumption eliminates third and higher order
uncertainties. Figure~\ref{fig:forecast} shows that an assumption of a
perfect model class is not well supported in our application to ensemble
weather forecasting, because during lead-times of 144 to~180 hours the
entire forecast ensemble fails to represent what actually happened. This
implies the assumption that the ensemble is random selection form the
distribution of possibilities is false. Nonetheless, table~\ref{tab:odds}
shows that our odds calculation is sufficiently conservative at these lead
times to cope fairly well with this level of uncertainty; pay-outs of
around 50 units compared to 200 units for capped probabilities. The goal,
however, was that the pay-outs should be zero on average. It would be
optimistic to suggest this were the case for lead times of 5 days or more.
Improvement of these odds forecasts is certainly possible. Either by
more careful consideration of the selection of the prior, or by abandoning
the simple frequency and Gaussian models for a more
sophisticated model that takes into account the conditional aspects of the
weather. In analogy to the urn game, a more sophisticated model means
looking more closely at the mixing process.

We have argued that odds forecasting has uses in investment and loss
mitigation, we claim also that it can be used for model assessment. The
results of table~\ref{tab:odds} suggest that the Gaussian model is quite
successful out to lead times of 5 days. Kernel-density based models may
well do better, extending to longer lead times, or having average pay-outs
closer to zero. Comparison with an ultimate-challenge client as we do
provides a diagnostic of the models performance, and furthermore, if a
model achieves a zero average pay-out, then excess of the odds provide
indications of the model's higher order uncertainty.

%%%%%%%%%%%%%%%%%%%%%%%%%%%%%%%%%%%%%%%%%%%%%%%%%%%%%%%%%%%%%%%%%%%%%%%%%%%
\bibliographystyle{plainnat}
\bibliography{refs}
%%%%%%%%%%%%%%%%%%%%%%%%%%%%%%%%%%%%%%%%%%%%%%%%%%%%%%%%%%%%%%%%%%%%%%%%%%%

\vspace{10cm} % forces first figure onto next page

\renewcommand{\baselinestretch}{1.0}

\clearpage

\begin{table}
  \centering
  \begin{tabular}{|c||c|c|c||c|c|c|}
    \hline
     Utility & \multicolumn{3}{c||}{Linear} &  \multicolumn{3}{|c|}{Logarithmic}\\ 
    \hline
    \hline
    Lead & 1. & 2. & $3'$. & 1. & 2. & $3'$. \\ 
    time & Freq. & Gaus. & Capped & Freq. & Gaus. & Capped \\
    (days) & odds  & odds  & Prob. & odds  & odds  & Prob.  \\
    \hline
1 &      7.3 &     47.6 &     36.2 &    -1.19 &    -4.03 &     4.94\\
2 &     65.5 &    -11.7 &      160 &     8.97 &       -2 &     12.5\\
3 &      117 &     14.9 &      233 &     20.3 &     5.21 &     21.6\\
4 &     78.9 &     3.75 &      175 &     16.6 &     5.28 &     20.8\\
5 &     70.9 &     68.2 &      186 &     14.1 &     4.54 &     18.8\\
6 &     73.7 &     7.71 &      166 &     15.6 &     6.29 &       20\\
7 &     80.5 &     82.6 &      196 &     20.5 &     13.1 &     26.2\\
8 &      105 &     24.3 &      238 &       21 &     9.89 &     27.6\\
9 &      124 &     50.1 &      283 &     28.5 &       18 &     34.4\\
10 &      147 &      130 &      325 &     33.5 &     22.9 &     38.1\\
    \hline
  \end{tabular}
  \caption{%{}
    Total pay-out of forecasters on bets taken from the 
    ultimate-challenge client on the event that the minimum 
    temperature at London Heathrow in the 24 hours period from 18:00UTC is
    more than 3 degrees centigrade below the control 
    temperature when the forecasts are issued. 
    Total computed over first half of 2005. In the case of logarithmic utility 
    the pay-out is expressed as the logarithm base 10.
  }
  \label{tab:odds}
\end{table}

\end{document}